\documentclass[11pt,leqno]{article}

\usepackage{amsmath,amsfonts,amscd,amssymb,theorem}

\long\def\comment#1\endcomment{}

\comment
\endcomment

\makeatletter
\begingroup
\gdef\th@dotted{\normalfont\itshape
  \def\@begintheorem##1##2{%
        \item[\hskip\labelsep \theorem@headerfont ##1\ ##2.]}%
\def\@opargbegintheorem##1##2##3{%
   \item[\hskip\labelsep \theorem@headerfont ##1\ ##2\ (##3).]}}
\endgroup
\makeatother

\theoremstyle{dotted}

\newtheorem{theorem}{Theorem}[section]
\newtheorem{lemma}[theorem]{Lemma}

\newtheorem{prop}[theorem]{Proposition}

\makeatletter
\begingroup
\gdef\th@upshape{\normalfont
  \def\@begintheorem##1##2{%
        \item[\hskip\labelsep \theorem@headerfont ##1\ ##2.]}%
\def\@opargbegintheorem##1##2##3{%
   \item[\hskip\labelsep \theorem@headerfont ##1\ ##2\ (##3).]}}
\endgroup
\makeatother

\theoremstyle{upshape}

\newtheorem{defn}[theorem]{Definition}
\newtheorem{remark}[theorem]{Remark}
\newtheorem{exa}[theorem]{Example}

\makeatletter
\renewcommand{\subsection}{\@startsection{subsection}{2}{0pt}{-3ex
plus -1ex minus -0.2ex}{-2mm plus -0pt minus
-2pt}{\normalfont\bfseries}} 
\renewcommand{\subsubsection}{\@startsection{subsubsection}{3}{0pt}{-3ex
plus -1ex minus -0.2ex}{-2mm plus -0pt minus
-2pt}{\normalfont\bfseries}} 
\makeatother

\makeatletter
\@addtoreset{equation}{section}
\makeatother

\newcommand{\cntrct}                
{\hspace{2pt}\raisebox{1pt}{\text{$\lrcorner$}}\hspace{2pt}}

\newcommand{\proof}[1][Proof.]{\smallskip\noindent{\em #1}}
\def\endproof{\hfill\ensuremath{\square}\par\medskip}

\def\eqref#1{\thetag{\ref{#1}}}

\let\latexref=\ref
\def\ref#1{{\normalfont{\latexref{#1}}}}

\newcommand{\wt}{\widetilde}

\newcommand{\idot}{{\:\raisebox{1pt}{\text{\circle*{1.5}}}}}
\newcommand{\hdot}{{\:\raisebox{3pt}{\text{\circle*{1.5}}}}}
\newcommand{\I}{\mathcal{I}}
\newcommand{\D}{{\cal D}}
\newcommand{\DF}{{\cal D}F}

\newcommand{\Z}{\mathbb{Z}}

\newcommand{\Coh}{\operatorname{Coh}}
\newcommand{\Fl}{\operatorname{Fl}}
\newcommand{\Bl}{\operatorname{Bl}}

\newcommand{\Proj}{\operatorname{Proj}}
\newcommand{\Spec}{\operatorname{Spec}}
\newcommand{\A}{\mathcal{A}}
\newcommand{\E}{\mathcal{E}}
\newcommand{\F}{\mathcal{F}}
\newcommand{\calo}{\mathcal{O}}

\newcommand{\amod}{{\text{\sf -mod}}}
\newcommand{\fmod}{{\text{\sf -filt}}}

\newcommand{\Id}{\operatorname{\sf Id}}
\newcommand{\id}{\operatorname{\sf id}}
\newcommand{\gr}{\operatorname{\sf gr}}

\title{Refined blowups}

\author{D. Kaledin\thanks{Partially supported by the Dynasty
    Foundation award and by a subsidy granted to the HSE by the
    Government of the Russian Federation for the implementation of
    the Global Competitiveness Program.}\ \ and 
  A. Kuznetsov\thanks{Partially supported by RFBR grants
    14-01-00416, 15-01-02164, 15-51-50045, NSh-2998.2014.1, and by a
    subsidy granted to the HSE by the Government of the Russian
    Federation for the implementation of the Global Competitiveness
    Program.}}

\begin{document}

\maketitle

\begin{abstract}
Categorical resolution of singularities has been constructed in
\cite{ku-lu}. It proceeds by alternating two steps of seemingly
different nature. We show how to use the formalism of filtered derived
categories to combine the two steps into one. This results in a certain
rather natural categorical refinement of the usual blowup of an algebraic
variety in a closed subscheme.
\end{abstract}

\section*{Introduction.}

Non-commutative algebraic geometry, at least in its modern
homological version, studies triangulated categories with some
enhancement (for example, DG enhancement as in \cite{kel}). Passing
from commutative to non-commutative world, then, consist of
associating such a category to any algebraic variety $X$. The usual
choice is to consider the bounded derived category $\D^b(X)$ of
coherent sheaves on $X$, and take its full subcategory $\D^{pf}(X)
\subset \D^b(X)$ formed by perfect complexes. While $\D^b(X)$
contains all the geometric objects such as the structure sheaves of
subvarieties $Z \subset X$, $\D^{pf}(X)$ better reflects the global
geometry of $X$. In particular, there are notions of homological
smoothness and homological properness for enhanced triangulated
categories, and $X$ is smooth resp. proper if and only if so is
$\D^{pf}(X)$. Moreover, if $X$ is smooth, $\D^{pf}(X)$ coincides
with the whole $\D^b(X)$.

Homologically smooth and proper categories are called saturated. It
is these categories that one usually studies in non-commutative
geometry, and it is only to these categories that many important
results apply.

If $X$ is singular, then $\D^{pf}(X)$ and $\D^b(X)$ are
different. If $X$ is proper, then $\D^{pf}(X)$ is homologically
proper but it is not homologically smooth. On the other hand,
$\D^b(X)$ is always homologically smooth by a recent result of Lunts
\cite{lu} using a surprising and deep result of Rouqier \cite{rou};
however, even if $X$ is proper, $\D^b(X)$ is not homologically
proper. Thus working with singular varieties presents a problem.

In the usual algebraic geometry, one handles problems of this type
by resolving the singularities of an algebraic variety $X$, that is,
constructing a smooth variety $Y$ equipped with a proper birational
map $Y \to X$. Categorical resolution of singularities introduced
recently in \cite{ku-lu} is a non-commutative counterpart and
refinement of this procedure. The goal is to embed $\D^{pf}(X)$
fully and faithfully into a category $\wt{\D}$ that behaves better
--- it is homologically smooth, homologically proper if $X$ is
proper, and can be represented as an iterated extension of
categories of the form $\D^b(Z)$ for smooth algebraic varieties
$Z$. Such a categorical resolution has been constructed and studied
in \cite{ku-lu} for all separated schemes of finite type over a
field of characteristic $0$. Taking a resolution of singularities
$f:Y \to X$ is the first step of the construction. However, the
pullback functor $L^\hdot f^*:\D^{pf}(X) \to \D^b(Y)$ is fully
faithful only if $X$ has rational singularities, so that in the
general case, much more work is required to achieve the goal.

Let us note that one of the main points of \cite{ku-lu} is the
interplay between commutative and non-\-com\-mu\-ta\-tive
geometry. From the purely non-commutative perspective, aside from
the fact that the resolution $\wt{\D}$ is homologically proper if
$X$ is proper, probably the main application of the construction is
a recent result of Efimov \cite{efi}. He uses the explicit form of
the construction of \cite{ku-lu} to refine Lunts's Theorem of
\cite{lu} and prove that $\D^b(X)$ is homologically of finite type
(``homologically finitely presented'' in the language of
\cite{to-va}). The fact that the homologically smooth category
$\wt{\D}$ is an iterated extension of categories of the form
$\D^b(Z)$, $Z$ smooth, plays no role --- it has no non-commutative
interpretation.

However, it would be presumptious to expect that the current general
non-commutative perspective is the final word in the subject: there
might well be some subtle ways in which categories coming from
smooth algebraic varieties are substantially different from general
homologically smooth categories.

In fact, if anything, our present knowledge argues exactly for this
--- there are deep purely non-commutative results that have only
been proved for saturated categories coming from algebraic
varieties. One example is the non-commutative Hodge-to-de Rham
degeneration theorem that has been proved in \cite{kal} under a
certain technical assumption, and we do not have a proof in the
completely general case. Another much more striking example is the
Blanc-To\"en Conjecture of \cite{kamo}: for categories of
commutative origin, it has been completely proved in \cite{blanc},
while the general case is wide open.

Thus at least at the moment, the conservative approach taken in
\cite{ku-lu} is obviously justified: knowing that a category is an
iterated extension of categories of commutative origin is quite
useful, both from the practical and from the conceptual point of
view.

\medskip

The point of the present paper is to show that the construction of
\cite{ku-lu} is also very natural and simple, in fact even more
natural than the usual commutative blowup.

\medskip

Namely, as given in \cite{ku-lu}, the categorical resulution
procedure works by alternating two main steps of completely
different nature. In this paper, we restrict our attention to one
sequence of step 1 followed by step 2, and observe that the two
steps can be combined into a single procedure.

More precisely, recall that the usual blowup $Y$ of an algebraic
variety $X$ in an ideal sheaf $\I$ is expressed as $Y = \Proj
\A_\idot$, where
$$
\A_\idot = \bigoplus_{n \geq 0}\I^n
$$
is a graded algebra sheaf on $X$. However, the sheaf $\A_\idot$ has
more structure --- it is the Rees algebra of a filtered algebra sheaf
(namely, the structure sheaf $\calo_X$ filtered by the powers $\I^n
\subset \calo_X$ of the ideal sheaf $\I$). We first show that the
usual Serre theorem expressing $\D^b(Y)$ in terms of graded sheaves
of modules over $\A_\idot$ works without any changes in the filtered
context, and then observe that after an elementary modification, the
filtered construction gives a ``partial categorical resolution'' ---
namely, it embeds $\D^{pf}(X)$ fully faithfully into an iterated
extension of a copy $\D^b(Y)$ and several copies of $\D^b(Z)$, where
$Z = \Spec(\calo_X/\I)$ is the center of the blowup.

To obtain the whole categorical resolution of singularities by this
construction, we need to know how to iterate it. It can probably be
done by considering sheaves equipped with several filtrations, but
at present, we have not seriously pursued this. Thus we only did a
very simple thing, and we try to present it as simply as possible
--- Section~\ref{1} contains the precise statements of the results,
and Section~\ref{2} is taken up with the proofs. It turns out that a
considerable simplification can be achieved by using the machinery
of exact categories and filtered derived categories, so we adopt
this approach. The reader who is uncomfortable with this machinery
can find an alternative formulation in terms of graded modules over
the Rees algebra in the beginning of Section~\ref{2}. Throughout the
paper, we restrict out attention to $\D^b(-)$ and completely ignore
the enhancements and the ``big'' unbounded derived categories of
quasicoherent sheaves. Taking care of either of these two points is
a straightforward exercise that does not seem to require new ideas.

\section{Statements.}\label{1}

Fix once and for all a Noetherian scheme $X$. Let $\Coh(X)$ be the
abelian category of coherent sheaves of $\calo_X$-modules, and let
$\D^b(X)$ be its bounded derived category. Let $\D^{pf}(X) \subset
\D^b(X)$ be the full subcategory spanned by perfect complexes. Let
$\Fl(X) \subset \Coh(X)$ be the full subcategory spanned by coherent
sheaves $E \in \Coh(X)$ that are flat over $\calo_X$ (or
equivalently, locally free). The category $\Fl(X)$ is an exact
category in the sense of Quillen.

In general, $\D^{pf}(X)$ is not the derived category of any abelian
category. However, assume from now on for simplicity that $\Coh(X)$
has enough flat coherent sheaves --- in other words, for any $\E \in
\Coh(X)$ there exists a surjective map $\E' \to \E$ with $\E' \in
\Fl(X)$ (for example, it suffices to assume that $X$ is
quasiprojective over an affine scheme). Then $\D^{pf}(X)$ is the
bounded derived category of the exact category $\Fl(X)$, as in
\cite{kel0}.

Assume given a scheme $Y$ equipped with a projective map $f:Y \to
X$. By definition, we have $Y = \Proj \A_\idot$, where $\A_\idot$ is
a sheaf of graded algebras on $X$ given by
$$
\A_n = f_*\calo_Y(n), \qquad n \geq 0,
$$
with $\calo_Y(n)$ being the powers of a relatively very ample line
bundle $\calo_Y(1)$ on~$Y$. For every $n$, $\A_n$ is a coherent sheaf
of $\calo_X$-modules, and $\A_\idot$ is generated by $\A_0$ and
$\A_1$ as an algebra sheaf.

Denote by $\A_\idot\amod$ the category of sheaves
$$
\E_\idot = \bigoplus_{i \in \Z}\E_i
$$
of graded $\A_\idot$-modules that are locally finitely generated
over $\A_\idot$ (in the paper, we will never consider non-graded
modules over a graded algebra, so that the notation is
unambiguous). Then the components $\E_n$ of every object $\E_\idot$
of the category $\A_\idot\amod$ are coherent sheaves of
$\calo_X$-modules, and the action map $\A_1 \otimes \E_i \to
\E_{i+1}$ is surjective for $i \gg 0$. The category $\A_\idot\amod$
is abelian. Denote by $\D^b(\A_\idot)$ its bounded derived
category. Let $\A_\idot\amod_{tors} \subset \A_\idot\amod$ be the
full abelian subcategory spanned by objects $\E_\idot$ such that
$\E_i = 0$ for $i \gg 0$, and let $\D^b_{tors}(\A_\idot)$ be its
bounded derived category. Then the subcategory $\A_\idot\amod_{tors}
\subset \A_\idot\amod$ is thick, the natural embedding
$\D^b_{tors}(\A_\idot) \subset \D^b(\A_\idot)$ is fully faithful,
and the famous theorem of Serre provides equivalences of categories
$$
\A_\idot\amod/\A_\idot\amod_{tors} \cong \Coh(Y), \qquad
\D^b(\A_\idot)/\D^b_{tors}(\A_\idot) \cong \D^b(Y),
$$
where $\Coh(Y)$ is the abelian category of coherent sheaves of
$\calo_Y$-modules, and $\D^b(Y)$ is its bounded derived
category\footnote{The abelian case is \cite{serre}. The derived case
  seems to be folklore, but see \cite[Lemma
    14, 15]{orlov}.}. Moreover, let us denote by
$\overline{f}^*:\A_\idot\amod \to \Coh(Y)$ the natural projection;
then it has a right-adjoint $\overline{f}_*:\Coh(Y) \to
\A_\idot\amod$ given by
\begin{equation}\label{wt.f}
(\overline{f}_*\E)_n = f_*\E(n), \qquad n \geq 0,
\end{equation}
and its derived functor $R^\hdot\overline{f}_*:\D^b(Y) \to
\D^b(\A_\idot)$ is right-adjoint to the projection
$\overline{f}^*:\D^b(\A_\idot) \to \D^b(Y)$. We have $\overline{f}^*
\circ R^\hdot \overline{f}_* \cong \id$, and we have a
semiorthogonal decomposition
$$
\D^b(\A_\idot) = \langle \D^b(Y),\D^b_{tors}(\A_\idot) \rangle
$$
in the sense of \cite{boka} (in other words, the projection
$\D^b(\A_\idot) \to \D^b(Y)$ induces an equivalence between the
category $\D^b(Y)$ and the right orthogonal
$\D^b_{tors}(\A_\idot)^\perp \subset \D^b(\A_\idot)$).

\medskip

Assume now that we are given an ideal sheaf $\I \subset \calo_X$,
and the algebra~$\A_\idot$ is actually given by its powers,
\begin{equation}\label{ku}
\A_n \cong \I^n, \qquad n \geq 0,
\end{equation}
so that $Y$ is the blowup of $X$ in the closed subscheme $Z \subset
X$ defined by the ideal sheaf $\I$. Here are two typical examples of
the situation.

\begin{exa}
If $X$ is quasiprojective and normal, and the projective map $f:Y
\to X$ is birational, then we can always find an ideal sheaf $\I
\subset \calo_X$ such that \eqref{ku} holds, so that $Y$ is
automatically a blowup. Indeed, we have $\A_0 = f_*\calo_Y \cong
\calo_X$ by Zariski Connectedness Theorem. Moreover, we can twist
the dual $\calo_Y(-1)$ to a relative very ample line bundle
$\calo_Y(1)$ by a pullback of a sufficiently ample line bundle on
$X$ so that it acquires a global section --- in other word, we have
an injective map $u:\calo_Y(1) \to \calo_Y$. Then $f_*(u):\A_1 \to
\A_0 = \calo_X$ identifies $\A_1$ with an ideal sheaf $\I \subset
\calo_X$, and for any $n \geq 2$, $f_*(u^n):\A_n \to \A_0$
identifies $\A_n$ with the $n$-th power $\I^n \subset \calo_X$.
\end{exa}

\begin{exa}\label{nilp.exa}
Conversely, assume that the ideal sheaf $\I \subset X$ is
nilpotent. Then the blowup $Y$ is perfectly well defined but empty.
\end{exa}

\begin{defn}
A {\em $Z$-filtered sheaf} $\langle \E,F^\hdot \rangle$ is a
coherent sheaf $\E$ of $\calo_X$-modules equipped with a decreasing
filtration by sheaves of $\calo_X$-submodules $F^n\E$, $n \geq 0$,
such that $\I \cdot F^n\E \subset F^{n+1}\E$ for any $n \geq 0$.
\end{defn}

For any $Z$-filtered sheaf $\langle \E,F^\hdot \rangle$ and any $n,m
\geq 0$, we then have $\I^m \cdot F^n\E \subset F^{n+m}\E$, and the
direct sum
\begin{equation}\label{rees}
\E_\idot = \bigoplus_{n \geq 0}F^n\E
\end{equation}
is a graded module over $\A_\idot = \bigoplus \I^n$. Denote by
$\A\fmod$ the category of $Z$-filtered sheaves $\langle \E,F^\hdot
\rangle$ such that $\E_\idot$ is locally finitely generated as a
sheaf of $\A_\idot$-modules (here $\A$ stands for the filtered
algebra sheaf $\langle \calo_X,\I^\hdot \rangle$, that is, $\calo_X$
equipped with the $\I$-adic filtration). Then $\A\fmod$ is an exact
category, and we can consider its derived category
$\DF^b(\A)$. Explicitly, say that a map $f$ of complexes of objects
in $\A\fmod$ is a {\em filtered quasiisomorphism} if the associated
graded map $\gr^\hdot(f)$ is a quasiisomorphism; then $\DF^b(\A)$
can be obtained by inverting filtered quasiisomorphisms in the
category of bounded complexes in $\A\fmod$, as in \cite{BBD}. There
is also a description in terms of graded modules over the so-called
{\em extended Rees algebra} that we recall in the beginning of
Section~\ref{2}.

\comment
\begin{remark}
By definition, $\A\fmod$ is the full subcategory in the
exact category of filtered coherent sheaves on $X$. Note, however,
that this subcategory is in general not closed under
extensions. Therefore the natural functor from $\DF^b(\A)$ to
the filtered version $\DF^b(X)$ of the derived category $\D^b(X)$ is
not fully faithful.
\end{remark}
\endcomment

For any locally free coherent sheaf $\E$ of $\calo_X$-modules,
equipping $\E$ with the filtration $F^n\E = \I^n \cdot \E \subset
\E$ turns it into a $Z$-filtered sheaf. The corresponding graded
module $\E_\idot$ is finitely generated (in fact, it is generated by
its degree-$0$ component), so that we obtain a natural functor
\begin{equation}\label{pb.0}
\Fl(X) \to \A\fmod.
\end{equation}
This functor is obviously exact, so that it induces a functor
\begin{equation}\label{pb}
\rho:\D^{pf}(X) = \D^b(\Fl(X)) \to \DF^b(\A).
\end{equation}

\begin{lemma}\label{rho.ff}
The functor $\rho$ of \eqref{pb} is fully faithful.
\end{lemma}

\proof{} Let $\A\fmod_{fl} \subset \A\fmod$ be the full subcategory
spanned by filtered sheaves $\langle \E,F^\hdot \rangle$ with $\E
\in \Fl(X)$. This subcategory is closed under extensions, so that it
is an exact subcategory, and its derived category $\DF^b_{fl}(\A)$
is a full subcategory in $\DF^b(\A)$. The embedding \eqref{pb.0}
obviously sends $\Fl(X)$ into $\A\fmod_{fl} \subset \A\fmod$, and
the forgetful functor
\begin{equation}\label{pb.2}
\A\fmod_{fl} \to \Fl(X)
\end{equation}
sending $\langle \E,F^\hdot \rangle$ to $\E$ is right-adjoint to the
embedding \eqref{pb.0}. Since the composition of \eqref{pb.0} and
\eqref{pb.2} is the identity functor, \eqref{pb.0} is fully
faithful. Moreover, \eqref{pb.2} is exact, thus induces a functor
$\DF^b_{fl}(\A) \to \D^{pf}(X)$ right-adjoint to \eqref{pb.0}, and
its composition with the functor $\rho$ of \eqref{pb} is again
isomorphic to the identity functor, so that $\rho$ is fully
faithful.
\endproof

We can now formulate our results. Denote by $\DF^b_{tors}(\A)
\subset \DF^b(\A)$ the full subcategory spanned by complexes of
modules $\E$ such that $F^n\E = 0$ for $n \gg 0$. Then sending
$\langle \E,F^\hdot \rangle$ to $\E_\idot$ of \eqref{rees} gives a
functor $\A\fmod \to \A_\idot\amod$. This functor extends to a
functor
\begin{equation}\label{rees.fu}
\DF^b(\A) \to \D^b(\A_\idot),
\end{equation}
sending $\DF^b_{tors}(\A)$ into $\D^b_{tors}(\A_\idot)$.

\begin{lemma}\label{1.le}
The functor \eqref{rees.fu} induces an equivalence of categories
$$
\DF^b(\A)/\DF^b_{tors}(\A) \cong
\D^b(\A_\idot)/\D^b_{tors}(\A_\idot) \cong \D^b(Y),
$$
and we have a semiorthogonal decomposition
$$
\DF^b(\A) = \langle \D^b(Y),\DF^b_{tors}(\A) \rangle.
$$
\end{lemma}

Moreover, say that a $Z$-filtered module $\langle \E,F^\hdot
\rangle$ is {\em $n$-stable} for some $n \geq 0$ if $F^n\E = \E$,
let $\DF^b(\A)_n \subset \DF^b(\A)$ be the full
subcategory spanned by complexes of $n$-stable sheaves, and let
$$
\DF^b_{tors}(\A)_n = \DF^b_{tors}(\A) \cap
\DF^b(\A)_n.
$$

\begin{lemma}\label{2.le}
For any $n \geq 0$, we have a natural equivalence of categories
$$
\DF^b_{tors}(\A)_n/\DF^b_{tors}(\A)_{n+1} \cong \D^b(Z),
$$
and the projection $\DF^b_{tors}(\A)_n \to \D^b(Z)$ has a
right-adjoint, so that we have a semiorthogonal decomposition
$$
\DF^b_{tors}(\A)_n = \langle
\D^b(Z),\DF^b_{tors}(\A)_{n+1} \rangle.
$$
\end{lemma}

\begin{defn}
For any integer $n \geq 0$, the {\em $n$-refined blowup} of $X$ in
$Z \subset X$ is the quotient category
$$
\Bl_n(X,Z) = \DF^b(\A)/\DF^b_{tors}(\A)_n.
$$
\end{defn}

Note that Lemma~\ref{1.le} and Lemma~\ref{2.le} imply that we have a
semiorthogonal decomposition
$$
\DF^b(\A) = \langle \Bl_n(X,Z),\DF^b_{tors}(\A)_n\rangle.
$$

\begin{prop}\label{prop}
For $n \gg 0$, the functor
$$
\rho_n:\D^{pf}(X) \to \Bl_n(X,Z) = \DF^b(\A)/\DF^b_{tors}(\A)_n
$$
obtained by composing the fully faithful embedding $\rho$ of
\eqref{pb} with the natural projection is fully faithful.
\end{prop}

We note that by Lemma~\ref{1.le}, the $0$-refined blowup
$\Bl_0(X,Z)$ is just the category $\D^b(Y)$. The functor $\rho_0$
then coincides with the natural pullback functor
\begin{equation}\label{f.st}
f^*:\D^{pf}(X) \to \D^b(Y).
\end{equation}
If $X$ has rational singularities, then already this functor is
fully faithful, and we have a categorical resolution even without
further refinements.

In the general case, however, we do need to take some $n \geq 1$
(see Remark~\ref{bound.rem} for an effective lower bound on $n$).
By Lemma~\ref{2.le}, for any such $n$, the $n$-refined blowup has a
natural descreasing filtration by full subcategories whose top
quotient is $\D^b(Y)$, and the other quotients are equivalent to
$\D^b(Z)$. Moreover, we in fact have a semiorthogonal decomposition
\begin{equation}\label{Y.ZZ}
\Bl_n(X,Z) = \langle \D^b(Y), \D^b(Z), \dots, \D^b(Z) \rangle,
\end{equation}
with $n$ copies of $\D^b(Z)$ in the right-hand
side. Proposition~\ref{prop} then shows that for a large enough $n$,
the $n$-refined blowup gives a ``partial categorical resolution'' of
$X$ in an appropriate sense. Note that in the situation of
Example~\ref{nilp.exa}, $\D^b(Y)$ is empty. However,
Proposition~\ref{prop} still works; in this case, it coincides with
the Auslander construction that forms Step 2 in the categorical
resolution procedure of \cite{ku-lu}.

\begin{remark}
One might be tempted to avoid the need for Proposition~\ref{prop} by
saying that $\DF^b(\A)$ itself is a categorical resolution of
$\D^{pf}(X)$. Indeed, the functor $\rho:\D^{pf}(X) \to \DF^b(\A)$ is
fully faithful by Lemma~\ref{rho.ff}, and it can be shown that
$\DF^b(\A)$ admits a version of the semiorthogonal decomposition
\eqref{Y.ZZ} with infinitely many copies of $\D^b(Z)$. However, such
``infinite'' semiorthogonal decompositions seem to be useless; in
practice, it is essential that we only have a finite number of
terms.
\end{remark}

\section{Proofs.}\label{2}

As in Section~\ref{1}, fix a Noetherian scheme $X$ with an ideal
sheaf $\I \subset \calo_X$, assume that $\Coh(X)$ has enough locally
free sheaves, let $Z \subset X$ be the closed subscheme defined by
$\I$, let
$$
\A_\idot = \bigoplus_{n \geq 0}\I^n,
$$
and let $Y = \Proj\A_\idot$. It is convenient to consider the {\em
  extended Rees algebra} $\wt{\A}_\idot$ associated to the $\I$-adic
filtration on $\calo_X$. This is a graded subalgebra $\wt{\A}_\idot
\subset \calo_X[u,u^{-1}]$ in the algebra $\calo_X[u,u^{-1}]$ of
Laurent polynomials in one variable $u$ of degree $-1$, and it is
given by
$$
\wt{\A}_n = 
\begin{cases}
\I^n, &\quad n \geq 1,\\
\calo_X, &\quad n \leq 0.
\end{cases}
$$
Alternatively, $\wt{\A}_\idot$ is the quotient of the polynomial
algebra $\A_\idot[u]$ by the relations $a \cdot u = u(a)$ for any $a
\in \A_1$, where $u:\A_1 = \I \to \A_0 = \calo_X$ is the natural
embedding. Yet another description is
$$
\wt{\A}_\idot = \A_\idot \oplus u\A_0[u],
$$
with the product given again by $a \cdot u = u(a)$. In particular,
we have an algebra embedding $\A_\idot \to \wt{\A}_\idot$ that
identifies $\A_\idot$ with the subalgebra of elements in
$\wt{\A}_\idot$ of non-negative degree.

The category $\wt{\A}_\idot\amod$ of graded sheaves of
$\wt{\A}_\idot$-modules locally finitely generated over
$\wt{\A}_\idot$ is then an abelian category, with the bounded
derived category $\D^b(\wt{\A}_\idot)$. Say that $\E_\idot \in
\wt{\A}_\idot\amod$ is {\em $n$-stable} for some integer $n$ if
$u:\E_i \to \E_{i-1}$ is an isomorphism for $i \leq n$, let
$\wt{\A}_\idot\amod_n \subset \wt{\A}_\idot\amod$ be the full
subcategory spanned by $n$-stable modules, and let
$\D^b_n(\wt{\A_\idot}) \subset \D^b(\wt{\A}_\idot)$ be its bounded
derived category. Then for any $Z$-filtered sheaf $\langle
\E,F^\hdot \rangle \in \A\fmod$, the graded sheaf $\E_\idot$ of
\eqref{rees} extends naturally to a $0$-stable graded module
$\wt{\E}_\idot \in \wt{\A}_\idot\amod$ given by
$$
\wt{\E}_n =
\begin{cases}
F^n\E, &\quad n \geq 1,\\
\E, &\quad n \leq 0,
\end{cases}
$$
with $u:\wt{\E}_\idot \to \wt{\E}_{\idot-1}$ acting by the natural
embedding. This defines an exact functor $\A\fmod \to
\wt{\A}_\idot\amod$ that identifies $\A\fmod$ with the full
subcategory in $\wt{\A}_\idot\amod$ spanned by $0$-stable modules
$\E_\idot \in \wt{\A}_\idot\amod$ such that the map $u:\E_i \to
\E_{i-1}$ is injective for any $i$. On the level of derived
categories, the corresponding functor
\begin{equation}\label{ree}
\DF^b(\A) \to \D^b_0(\wt{\A}_\idot)
\end{equation}
is an equivalence of categories. For any $n \geq 1$, it sends
$n$-stable $Z$-filtered modules to $n$-stable
$\wt{\A}_\idot$-modules and identifies $\DF^b(\A)_n$ with
$\D^b_n(\wt{\A}_\idot)$. Under the equivalence \eqref{ree}, the
functor $\rho$ of \eqref{pb} sends $\E \in \Fl(X)$ to $\E
\otimes_{\calo_X} \wt{\A}_\idot$. The functor
$$
\tau:\D^b_0(\wt{\A}_\idot) \to \D^b(\A)
$$
corresponding to \eqref{rees.fu} is induced by the exact forgetful
functor
\begin{equation}\label{tau}
\tau:\wt{\A}_\idot\amod_0 \to \A_\idot\amod
\end{equation}
forgetting the negative degree components and the action of the
generator $u \in \wt{\A}_\idot$.

\proof[Proof of Lemma~\ref{1.le}.] Note that the functor
$\overline{f}_*:\Coh(Y) \to \A\amod$ of \eqref{wt.f} factors
through the projection $\tau$ of \eqref{tau}. Indeed, the natural
maps $u:\calo_Y(n+1) \to \calo_Y(n)$, $n \geq 0$ induce maps
\begin{equation}\label{new.u}
u:(\overline{f}_*\E)_{n+1} \to (\overline{f}_*\E)_n,
\end{equation}
for any $\E \in \Coh(Y)$, so that if we set
$$
\wt{f}_*(\E)_n = 
\begin{cases}
\overline{f}_*(\E)_n, &\quad n \geq 1,\\
\overline{f}_*(\E)_0, &\quad n \leq 0,
\end{cases}
$$
then the natural $\A_\idot$-action on $\overline{f}_*(\E)$ and the
maps \eqref{new.u} equip the graded sheaf $\wt{f}_*(\E)_\idot$ with
an action of the extended Rees algebra $\wt{\A}_\idot$. Thus we
obtain a functor
$$
\wt{f}_*:\Coh(Y) \to \wt{\A}_\idot\amod_0
$$
such that $\overline{f}_* \cong \tau \circ \wt{f}_*$, where $\tau$
is the forgetful functor \eqref{tau}.

Moreover, for any $\E_\idot \in \wt{\A}_\idot\amod_0$, the
adjunction map $\tau(\E_\idot) \to
\overline{f}_*\overline{f}^*\tau(\E_\idot)$ is compatible with the
maps \eqref{new.u}, thus comes from a map $\E_\idot \to
\wt{f}_*\wt{f}^*\E$, where we denote $\wt{f}^* = \overline{f}^*
\circ \tau$.

Passing to the derived functors, we obtain a functor
$$
R^\hdot\wt{f}_*:\D^b(Y) \to \D^b_0(\wt{\A}_\idot)
$$
such that $\tau \circ R^\hdot\wt{f}_* \cong R^\hdot\overline{f}^*$,
and a functorial map
\begin{equation}\label{adj}
a:\E_\idot \to R^\hdot\wt{f}_* \wt{f}^*\E_\idot
\end{equation}
for any $\E_\idot \in \D^b_0(\A_\idot)$. Note that we have
\begin{equation}\label{iso}
\wt{f}^* \circ R^\hdot\wt{f}_* = \overline{f}^* \circ \tau
\circ R^\hdot\wt{f}_* \cong \overline{f}^* \circ R^\hdot\overline{f}_*
\cong \Id,
\end{equation}
and $\tau(a)$ is the adjunction map between $\overline{f}^*$ and
$R^\hdot\overline{f}_*$. Then it is easy to check that \eqref{adj}
together with the isomorphism \eqref{iso} define an adjunction
between the functors $\wt{f}^*$ and $R^\hdot\wt{f}_*$, so that
\eqref{iso} implies that $R^\hdot\wt{f}^*$ is a fully faithful
embedding onto a left-admissible subcategory. To finish the proof, it
suffices to notice that the kernel of the functor
$$
\wt{f}^*:\DF^b(\A) \cong \D^b_0(\wt{\A}_\idot) \to \D^b(Y)
$$
is exactly the subcategory $\DF^b_{tors}(\A) \subset \DF^b(\A)$.
\endproof

\proof[Proof of Lemma~\ref{2.le}.] Denote by $i:Z \to X$ the
embedding. For any coherent sheaf $\E$ on $Z$, $i_*\E$ with the
filtration $F^ni_*\E = i_*\E$, $F^{n+1}i_*\E = 0$ is an $n$-stable
$Z$-filtered sheaf on $X$. This gives a fully faithful exact functor
$i_n:\Coh(Z) \to \A\fmod_n$ to the full subcategory $\A\fmod_n
\subset \A\fmod$ spanned by $n$-stable $Z$-filtered sheaves.

The functor $\gr^n_F:\A\fmod_n \to \Coh(Z)$ sending a $Z$-filtered
sheaf to its $n$-th associated graded quotient with respect to the
filtration is obviously left-adjoint to the functor $i_n$. The
functors $i_n$ and $\gr^n_F$ then extend to an adjoint pair of
functors between the categories of complexes of objects in the exact
categories $\A\fmod_n$, $\Coh(Z)$. Moreover, $i_n$ is exact, and the
functor $\gr^n_F$ is by definition also exact with respect to the
exact structure on $\A\fmod$, that is, it sends filtered
quasiisomorphisms between filtered complexes to
quasiisomorphisms. Thus $i_*$ and $\gr^n_F$ descend to an adjoint
pair of functors
$$
i_n:\D^b(Z) \to \DF^b(\A)_n, \qquad \gr^n_F:\DF^b(\A)_n \to
\D^b(Z),
$$
between derived categories, and the isomorphism $\gr^n_F \circ i_n
\cong \Id$ on the level of complexes induces an analogous
isomorphism on the level of derived categories. Therefore
$i_n:\D^b(Z) \to \DF^b(\A)_n$ is fully faithful.  It remains to
notice that $i_n$ sends $\D^b(Z)$ into $\DF^b_{tors}(\A)_n \subset
\DF^b(\A)_n$, and the orthogonal ${}^\perp i_n(\D^b(Z)) \subset
\DF^b_{tors}(\A)_n$ --- that is, the kernel of the functor $\gr^n_F$
--- is exactly $\DF^b_{tors}(\A)_{n+1}$.
\endproof

\begin{remark}
In particular, for $n=0$, the proof of Lemma~\ref{2.le} shows that
the functor
$$
i_0:\D^b(Z) \to \DF^b(\A)
$$
is fully faithful. Note that if we compose $i_0$ with the
tautological forgetful functor $\DF^b(\A) \to \D^b(X)$, $\langle
\E,F^\hdot \rangle \mapsto \E$, then the composition is the direct
image functor $i_*$ with respect to the embedding $i:Z \to X$, and
it certainly is not fully faithful on the level of derived
categories. Thus working with filtered derived categories is crucial
for the construction.
\end{remark}

Before we prove Proposition~\ref{prop}, let us make one
observation. Lemma~\ref{1.le} provides a natural quotient functor
$$
\wt{f}^*:\DF^b(\A) \to \DF^b(\A)/\DF^b_{tors}(\A)
\cong \D^b(Y),
$$
and as noted in Section~\ref{1}, the composition $\rho_0 = \wt{f}^*
\circ \rho:\D^{pf}(X) \to \D^b(Y)$ is isomorphic to the pullback
functor $f^*$ of \eqref{f.st}. More generally, for any object $\E \in
\Fl(X)$ and any filtered complex $\F \in \DF^b(\A)$, the tensor
product $\E \otimes_{\calo_X} \F$ is naturally a filtered complex,
and we have a natural functorial identification
\begin{equation}\label{pb.t}
\wt{f}^*(\E \otimes_{\calo_X} \F) \cong f^*\E \otimes_{\calo_Y}
\wt{f}^*\F.
\end{equation}
Lemma~\ref{1.le} also provides a functor
$R^\hdot\wt{f}_*:\D^b(Y) \to \DF^b(\A)$ right-adjoint to $\wt{f}^*$,
and for any $\E \in \Fl(X)$, $\F \in \Coh(Y)$, we have a natural map
\begin{equation}\label{adj.proj}
\begin{CD}
\wt{f}^*(\E \otimes_{\calo_X} R^\hdot\wt{f}_*\F) @>{a}>> f^*\E
\otimes_{\calo_Y} \wt{f}^*R^\hdot\wt{f}_*\F @>{\id \otimes b}>> f^*\E
\otimes_{\calo_Y} \F,
\end{CD}
\end{equation}
where $a$ is the isomorphism \eqref{pb.t}, and
$b:\wt{f}^*R^\hdot\wt{f}_*\F \to \F$ is the adjunction map. Then by
adjunction, \eqref{adj.proj} induces a map
\begin{equation}\label{proj}
\E \otimes_{\calo_X} R^\hdot\wt{f}_*\F \to R^\hdot\wt{f}_*(\wt{f}^*\E
\otimes_{\calo_Y} \F),
\end{equation}
and we observe that \eqref{proj} is an isomorphism. Indeed, it
suffices to check it separately on every filtered component, and
then it immediately follows from the usual projection formula and
\eqref{rees}.

\proof[Proof of Proposition~\ref{prop}.] By Lemma~\ref{1.le} and
Lemma~\ref{2.le}, the projection $\DF^b(\A) \to \Bl_n(X,Z)$
actually comes from a semiorthogonal decomposition. In particular,
it has a right-adjoint functor identifying $\Bl_n(X,Z)$ with the
right orthogonal
$$
\DF^b_{tors}(\A)_n^\perp \subset \DF^b(\A).
$$
Since the functor $\rho$ of \eqref{pb} is fully faithful, it
suffices to prove that for $n \gg 0$, it sends the whole
$\D^{pf}(X)$ into this orthogonal. 

Consider first the object $\calo_X \in \D^{pf}(X)$, and let
$\overline{\calo}_X \in \DF^b(\A)$ be the cone of the adjunction map
$\rho(\calo_X) \cong \A \to R^\hdot\wt{f}_*\wt{f}^*\A \cong
R^\hdot\wt{f}_*\calo_Y$. We then have a distinguished triangle
\begin{equation}\label{tria}
\begin{CD}
\A @>>> R^\hdot\wt{f}_*\calo_Y @>>> \overline{\calo_X} @>>> \A[1],
\end{CD}
\end{equation}
and by definition, $\overline{\calo}_X$ lies inside
$\DF^b_{tors}(\A) \subset \DF^b(\A)$, so that
$\gr^n_F\overline{\calo}_X = 0$ for $n \geq N$ for some fixed
constant $N \geq 0$.

Now for every $\E \in \Fl(X)$, the product $\E \otimes_{\calo_X}
\overline{\calo}_X$ also lies in the subcategory $\DF^b_{tors}(\A)
\subset \DF^b(\A)$, and we have
$$
\gr^n_F(\E \otimes_{\calo_X} \overline{\calo}_X) \cong \E \otimes_{\calo_X}
\gr^n_F\overline{\calo}_X = 0
$$
for $n \geq N$. By Lemma~\ref{2.le}, this implies
that $\E \otimes_{\calo_X} \overline{\calo}_X$ lies in
$$
\DF^b_{tors}(\A)_n^\perp \subset \DF^b_{tors}(\A)
\subset \DF^b(\A).
$$
On the other hand, by \eqref{proj}, the product $\E
\otimes_{\calo_X} R^\hdot\wt{f}_*\calo_Y$ lies in the image of the
fully faithful embedding $R^\hdot\wt{f}_*$. Therefore it is
orthogonal to the category $\DF^b_{tors}(\A)$ containing
$\DF^b_{tors}(\A)_n$, that is, also lies in
$\DF^b_{tors}(\A)_n^\perp$.

Tensoring $\E$ with the exact triangle \eqref{tria}, we deduce that
for any object $\E \in \Fl(X)$, hence also for any $\E \in
\D^{pf}(X)$, the object $\rho(\E) = \E \otimes_{\calo_X} \A$ lies in
$\DF^b_{tors}(\A)_n^\perp$ as soon as $n \geq N$.
\endproof

\begin{remark}\label{bound.rem}
As we see from the proof, the constant $n$ in Proposition~\ref{prop}
admits an effective lower bound: it is necessary and sufficient that
$$
R^\hdot\calo_Y(m) \cong \I^m
$$
for any $m \geq n$.
\end{remark}

We conclude with remarking that throughout the paper, the actual
geometry of the schemes $X$ and $Y$ has been used very little in
either the statements or the proofs. We in fact only need it to
control the size of things --- namely, to show that various
categories of finitely generated modules are abelian, and that the
adjoint functor $\overline{f}_*$ of \eqref{wt.f} exists on the level of
bounded derived categories. If we drop the finiteness and
boundedness assumptions, then everything works for arbitrary sheaves
of graded associative algebras, in fact even for non-commutative
ones (for Lemma~\ref{2.le}, the same proof goes through literally,
and for Lemma~\ref{1.le} and Proposition~\ref{prop}, the arguments
are very easy to generalize). One can also get by with weaker
conditions on the algebras such as coherence. We did not pursue this
since, as we have explained in the introduction, our main motivation
is in any case geometric, and at the moment, we do not have
interesting examples where larger generality is required.

\subsection*{Acknowledgements.} This paper was originally conceived
during the conference in Tokyo University celebrating the 60th
anniversary of Prof. Yu. Kawamata --- the second author was giving a
talk on \cite{ku-lu}, and the first author, while chairing the talk,
was trying to come up with a meaningful question. In the event,
there was so much interest in the talk that the question was not
required, so it was only asked later and in private. The question
was, is it possible to combine Step 1 and Step 2 of
\cite{ku-lu}. The answer is the present paper. We take this
opportunity to congratulate Prof. Kawamata once again, and to thank
the organizers of the conference, in particular Prof. K. Oguiso, for
inviting us both and for setting up this opportune moment.

\bigskip

\noindent
{\sc
Steklov Math Institute, Algebraic Geometry section\\
\mbox{}\hspace{30mm}and\\
Laboratory of Algebraic Geometry, NRU HSE}

\smallskip

(both institutions for both authors)

\medskip

\noindent
{\sc Center for Geometry and Physics, IBS,\\
Pohang, Rep. of Korea} (for D.K.)

\medskip

\noindent
{\em E-mail addresses\/}: {\tt kaledin@mi.ras.ru}, {\tt akuznet@mi.ras.ru}


\begin{thebibliography}{666}
{\footnotesize

\bibitem[BBD]{BBD} A. Beilinson, J. Bernstein, and P. Deligne, {\em
  Faisceaux Pervers}, Ast\'erisque {\bf 100}, Soc. Math. de France,
  1983.

\bibitem[Bl]{blanc} A. Blanc, {\em Invariants topologiques des
  Espaces non commutatifs}, arXiv:1307.6430.

\bibitem[BK]{boka} A. Bondal and M. Kapranov, {\em Representable
functors, Serre functors, and reconstructions}, (Russian) Izv. Akad.
Nauk SSSR Ser. Mat. {\bf 53} (1989), 1183--1205, 1337; translation
in Math. USSR-Izv. {\bf 35} (1990), 519--541.

\bibitem[E]{efi} A. Efimov, {\em Homotopy finiteness of some DG
  categories from algebraic geometry}, arxiv:1308.0135.

\bibitem[Ka1]{kal} D. Kaledin, {\em Non-commutative Hodge-to-de Rham
  degeneration via the method of Deligne-Illusie}, Pure
  Appl. Math. Q. {\bf 4} (2008), 785--875.

\bibitem[Ka2]{kamo} D. Kaledin, {\em Motivic structures in
  non-commutative geometry}, Proc. ICM 2010.

\bibitem[Ke1]{kel0} B. Keller, {\em Derived categories and universal
  problems}, Communications in Algebra {\bf 19} (1991), 699-747.

\bibitem[Ke2]{kel} B. Keller, {\em On differential graded
    categories}, in {\em  International Congress of Mathematicians},
  Vol. II, Eur. Math. Soc., Z\"urich, 2006; 151--190.

\bibitem[KL]{ku-lu} A. Kuznetsov and V. Lunts, {\em Categorical
  resolutions of irrational singularities}, arxiv:1212.6170.

\bibitem[L]{lu} V. Lunts, {\em Categorical resolution of
  singularities}, J. Algebra {\bf 323} (2010), 2977--3003.

\bibitem[O]{orlov} D. Orlov, {\em Derived categories of coherent
  sheaves and triangulated categories of singularities}, in Algebra, Arithmetic, and Geometry, 
  pp. 503-531. Birkhäuser Boston, 2009.
  
\bibitem[R]{rou} R. Rouquier, {\em Dimensions of triangulated categories}, 
Journal of K-theory: K-theory and its Applications to Algebra, Geometry, and Topology (1) {\bf 02} (2008), 193--256.
  
  
  
\bibitem[S]{serre} J.-P. Serre, {\em Faisceaux alg\'ebriques
  coherents}, Ann. of Math. {\bf 61} (1955), 197-278.

\bibitem[TV]{to-va} B. To\"en and M. Vaqui\'e, {\em Moduli of objects
    in dg-categories},  Ann. Sci. \'Ecole Norm. Sup. (4)  {\bf 40}
  (2007), 387--444.

}
\end{thebibliography}
\end{document}